# An image-incorporated immersed boundary method for diffusion equations

Andrew C. Chiang[a], Po-Nan Li[b,c], Soichi Wakatsuki[a,c]

[a] Department of Structural Biology, Stanford University, Stanford, CA 94305, USA

[b] Department of Electrical Engineering, Stanford University, Stanford, CA 94305, USA

[c] SLAC National Accelerator Laboratory, Menlo Park, CA 94025, USA

**Corresponding authors:** Andrew C. Chiang and Soichi Wakatsuki
Email: andrewchiang2011@gmail.com; soichi.wakatsuki@stanford.edu
Address: 318 Campus Drive West, Clark Center, Room W250, Stanford, CA 94305-5151, USA

**Abstract**

A novel sharp interface ghost-cell based immersed boundary method has been proposed and its parameters have been optimized against an analytical model in diffusion applications. The proposed embedded constrained moving least-squares (ECMLS) algorithm minimizes the error of the interpolated concentration at the image point of the ghost point by applying a moving least-squares method on all internal nodes, near the ghost image point, and the associated mirrored image points of these internal nodes through the corresponding boundary conditions. Using an analytical model as a reference, the ECMLS algorithm is compared to the constrained moving least-squares (CMLS) algorithm and the staircase model using various grid sizes, interpolation basis functions, weight functions, and the penalty parameter of the constraint. It is found that using ECMLS algorithm in the investigated diffusion application, the incomplete quartic basis function yields the best performance while the quadratic, cubic, and bicubic basis functions also produce results better than the staircase model. It is also found that the linear and bilinear basis functions cannot produce results better than the staircase model in diffusion applications. It is shown that the optimal radius of the region of internal nodes used for interpolation scales with the logarithm of the boundary radius of curvature. It is shown that for the diffusion application, the proposed ECMLS algorithm produces lower errors at the boundary with better numerical stability over a wider range of basis functions, weight functions, boundary radius of curvature, and the penalty parameter than the CMLS algorithm.

## 1. Introduction

The finite difference method is commonly used in numerically solving partial differential equations because of the ease in discretization and approximation of derivatives using algebraic equations of variable values at grid nodes [1] [2] [3]. However, when the shape of the boundary is irregular and does not conform to the grid lines, special procedures are needed to model the boundary conditions properly [3] [4] [5]. The simplest way to model a non-conforming boundary is to use a staircase model where the grid points outside of the boundary are removed and the boundary conditions are applied on the staircase grid lines [6]. One-dimensional boundary condition in each Cartesian coordinate is applied independently in the staircase model. Despite its inaccuracy, the staircase model is still used in many applications due to its simplicity and convenience [7] [8].

One way to enforce boundary conditions is to use the body-fitted grid methods where the governing equations and boundary conditions are transformed to body-fitted coordinate systems that conform to the boundaries [9] [10]. The grid-generation process of body-fitted structured or unstructured grids for complex geometries can be complicated and time-consuming because the mesh quality can impact error performance [11] [12]. Furthermore, for applications involving moving boundaries, transient remeshing procedures are needed which further increase the algorithmic and computational complexity of the body-fitted grid methods.

The immersed boundary method provides a framework to model partial differential equations in Cartesian grids and handle boundaries and interfaces also in the context of the Cartesian grid [13] [14] [15]. A major advantage of the immersed boundary method is that the grid generation does not need to conform to the geometry of the boundaries. However, since the boundary can cut through the mesh in arbitrary ways, it is critical to model the boundary properly to preserve the accuracy and conservation property of the numerical solver.

One type of the immersed boundary method is the so-called "continuous forcing" approach or diffused interface method where the immersed boundary is smeared by distributing singular Dirac delta forces over a thin layer of grid nodes in the vicinity of the boundary [13] [14]. In the diffused interface method, the boundary conditions are satisfied "gradually" over a localized region near the boundary instead of at the exact location of the boundary.

Another type of the immersed boundary method is the sharp interface method which includes the cut-cell method [16] [17] [18] and the ghost-cell method [19] [20]. The sharp interface method is preferred in the application where accurate solution near the boundary is desirable. In the cut-cell methods, the solid body is cut along the boundary to create a fraction of a cell. The finite volume is incorporated into the solver so that precise conservation of mass, momentum, and energy can be achieved. On the other hand, the ghost-cell methods use information at the internal nodes and the proximity of internal nodes to the boundary to interpolate or extrapolate the information at the ghost cell, which is outside of the boundary, but its information is needed for the calculation of the solution. The cut-cell methods are generally more accurate and produce less spurious fluctuations than the ghost-cell methods. However, the cut cells might be arbitrarily small and can lead to numerical stability issues if the time step is not reduced. The ghost-cell methods are preferred in many applications due to the simplicity in implementation and efficiency in computation because the complicated cell reshaping procedure is needed in cut-cell methods but not in ghost-cell methods.

Despite the simplicity in implementation, the most critical issue with the ghost-cell method is its ability to accurately interpolate or extrapolate the value at the ghost cell from the internal nodes and the boundary values. The interpolation method affects the number of points needed for the interpolation, the computation complexity to find the coefficients and to do interpolation, and the accuracy of interpolation. The simplest models use bilinear interpolation basis functions for two-dimensional (2-D) problems [21] [22] and trilinear basis functions for 3-D problems [23] [24]. To get better accuracy, the second-order quadratic basis function has also been used for interpolation [16] [19] [20]. The number of points needed for interpolation is equal to the rank of the interpolation basis functions. It is typically beneficial to use more data points to perform interpolation. The moving least-squares (MLS) technique has been used to minimize the sum of squared errors at all data points to find the best interpolation coefficients for trilinear basis functions [24] and quadratic basis functions [25] [26]. To further improve accuracy, higher-order polynomial basis functions are also used with MLS algorithm to perform interpolation [27] [28] [29]. When a large number of points over a wide area are used for interpolation, especially using the MLS algorithm, a smooth monotonic weight function is typically used to assign higher weights to points closer to the point of interpolation when the sum of squared errors is evaluated. Cosine weight functions [27] [28] [30], cubic spline weight functions [31], and inverse distance weight function [26] have been used as the weight function for MLS algorithm, while no weight functions are used in other implementations [24] [29]. The main goal of the ghost-cell method is using the values at internal nodes to accurately interpolate or extrapolate the value at the ghost cell which is outside of the boundary. This process can lead to inaccurate or numerically unstable solutions. The image point which is the ghost point mirrored from the boundary is introduced to address this issue [20] [24] [27] [30] [31] because the mirrored image point is inside the boundary, closer to the internal points used for interpolation. The value at the ghost point is related to the value at the image point through the boundary conditions at the projected point at the boundary.

Without using the boundary relationship between ghost point and image point, the MLS algorithm minimizes the total squared errors at the internal nodes used in interpolation. However, the MLS algorithm can be numerically unstable if the number of points used in interpolation is not sufficiently larger than the number of coefficients needed for interpolation [27]. The constrained moving least-squares (CMLS) algorithm [31] adds the boundary constraint between the ghost point and image point to the squared error minimization process to improve numerical stability. In this paper, an embedded constrained moving least-squares (ECMLS) algorithm is proposed to further improve numerical stability. In the ECMLS algorithm, all internal nodes used for interpolation are mirrored across the boundary to create more image points for interpolation. Each internal point and the associated image point are related through the boundary condition. The effectiveness and accuracy of the ECMLS algorithm are investigated in the paper using an analytical solution as a reference. The effects of various interpolation basis functions and weight functions are also investigated. The efficacy of the ECMLS algorithm is compared to the staircase model and the CMLS algorithm.

## 2. Two-dimensional Finite Difference Diffusion Model with an Immersed Boundary Method

Fick's second law of diffusion in two-dimensional homogeneous space can be described by the partial differential equation

$$\frac{\partial c}{\partial t}(x,y,t) = D\left[\frac{\partial^2 c}{\partial x^2}(x,y,t) + \frac{\partial^2 c}{\partial y^2}(x,y,t)\right], \tag{1}$$

where $c(x,y,t)$ is the concentration and $D$ is the diffusion coefficient. A Neumann boundary condition specifies the flux normal to the boundary

$$\psi_n(\Gamma,t) = -D\frac{\partial c}{\partial \hat{n}}(\Gamma,t) = \psi_{BP}(\Gamma), \tag{2}$$

where $\hat{n}$ is the direction normal to the boundary and $\psi_{BP}(\Gamma)$ is the flux normal to the boundary at the boundary $\Gamma$. The Neumann boundary condition (2) can be rewritten as

$$\frac{\partial c}{\partial \hat{n}}(\Gamma,t) = -\frac{1}{D}\psi_{BP}(\Gamma). \tag{3}$$

On the other hand, a Dirichlet boundary condition at boundary $\Gamma$ specifies the concentration at the boundary point

$$c(\Gamma,t) = c_{BP}(\Gamma), \tag{4}$$

where $c_{BP}(\Gamma)$ is the concentration at the boundary point. The Dirichlet boundary condition (4) implies that the second derivative of concentration in the normal direction $\hat{n}$ is zero at the boundary,

$$\frac{\partial^2 c}{\partial \hat{n}^2}(\Gamma,t) = \frac{1}{D}\frac{\partial c}{\partial t}(\Gamma,t) = \frac{1}{D}\frac{\partial}{\partial t}c_{BP} = 0. \tag{5}$$

The diffusion equation (1) can be modeled numerically using finite difference method with the forward time centered space (FTCS) representation [2]

$$\frac{c_{j,k}^{n+1} - c_{j,k}^n}{\Delta t} = D\left[\frac{c_{j+1,k}^n - 2c_{j,k}^n + c_{j-1,k}^n}{(\Delta x)^2} + \frac{c_{j,k+1}^n - 2c_{j,k}^n + c_{j,k-1}^n}{(\Delta y)^2}\right], \tag{6}$$

where $n$ is the time index, $\Delta t$ is the step size in time, $j$ is the space index in $x$ direction, $\Delta x$ is the step size in $x$ direction, $k$ is the space index in $y$ direction, $\Delta y$ is the step size in $y$ direction, and $c_{j,k}^n$ represent the discrete sample of concentration at $c(j\Delta x, k\Delta y, n\Delta t)$. The equation (6) can be rewritten as

$$c_{j,k}^{n+1} = c_{j,k}^n + \frac{D\Delta t}{(\Delta x)^2}\left(c_{j+1,k}^n - 2c_{j,k}^n + c_{j-1,k}^n\right) + \frac{D\Delta t}{(\Delta y)^2}\left(c_{j,k+1}^n - 2c_{j,k}^n + c_{j,k-1}^n\right). \tag{7}$$

Using equation (7), the samples at timestep $n$ can be used to evaluate $c_{j,k}^{n+1}$ at timestep $n+1$. The stability criterion is [2]

$$\frac{D\Delta t}{(\Delta x)^2} \leq \frac{1}{2}, \quad \frac{D\Delta t}{(\Delta y)^2} \leq \frac{1}{2}. \tag{8}$$

The nodes in Cartesian coordinates are partitioned into internal domain $\Omega_i$ and external domain $\Omega_e$ with the boundary $\Gamma$ separating the two domains (Fig. 1). Concentration values at the points in the internal domain are calculated using equation (7). The grid points in the internal domain can be classified as internal nodes (IN) and internal boundary nodes (IBN). An IN has all the points in the finite difference stencil (7) within the internal domain. An IBN has some points in the finite difference stencil (7) outside of the boundary. Those external points needed for computation are ghost points (GP). The concentration values at GPs are not updated based on equation (7). The other points in the external domain are external nodes (EN).

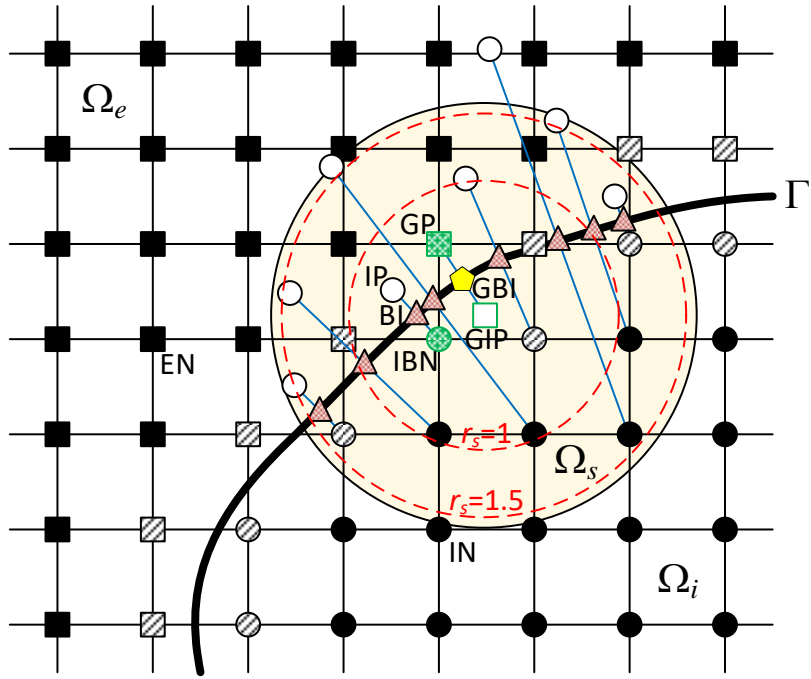

Fig. 1. A two-dimensional Cartesian grid with an immersed boundary; ● internal node (IN); ⊘ internal boundary node (IBN); ▨ ghost point (GP); ■ external node (EN); ▲ boundary intercept (BI); ○ image point (IP); ⊛ target IBN; ▦ target GP; ⬠ ghost boundary intercept (GBI); ☐ ghost image point (GIP). The normalized radius is given by $r_s = r/\sqrt{\Delta x^2 + \Delta y^2}$.

Consider the target IBN where an associated target GP is needed in equation (7). The target GP is normally projected onto the boundary at the ghost boundary intercept (GBI) and is mirrored from the GBI to the ghost image point (GIP). The value at the target GP is related to the value at the GIP though boundary conditions. Using the ghost-cell based sharp interface immersed boundary method, the value at the GIP is estimated using the points in the internal domain and the boundary conditions. With the MLS algorithm [27], the INs and IBNs near the GBI and inside the interpolation stencil domain $\Omega_s$ are used to interpolate the value at GIP by minimizing the weighted squared errors of the interpolating function at the INs and IBNs. Addressing the stability issue of the MLS algorithm, the constrained moving least-squares (CMLS) algorithm [31] adds the constraint of the boundary condition at the GBI into the error minimization equation in the MLS algorithm. In the currently proposed embedded constrained moving least-squares (ECMLS) algorithm, each of the INs and IBNs inside the interpolation stencil domain $\Omega_s$ associated with the target IBN is normally projected onto the boundary at a boundary intercept (BI) and is mirrored from the BI to an image point (IP). The values at the BI and IP are associated with the value at the IN or IBN through boundary conditions. The values at the IPs, BIs, and the INs and IBNs in $\Omega_s$ are used to minimize error of the interpolation function in order to estimate the values at GIP and target GP. The interpolation stencil domain $\Omega_s$ used in the ECMLS algorithm is centered at the GIP.

Using the Neumann boundary condition (3), the value $c_{IP}^n$ at the IP can be approximated by

$$c_{IP}^n \approx c_{IN}^n - \frac{2\delta}{D}\psi_{BI}, \tag{9}$$

where $\psi_{BI}$ is the flux at the BI, $\delta$ is the distance from the IN to the BI, and $c_{IN}^n$ is the concentration at the IN. Using the Dirichlet boundary condition (4) and (5), the value $c_{IP}^n$ at the IP can be approximated by

$$c_{IP}^n \approx -c_{IN}^n + 2c_{BI}, \tag{10}$$

where $c_{BI}$ is the concentration at the BI. Equations (9) and (10) can be consolidated into a common form

$$c_{IP}^n \approx \eta c_{IN}^n + \gamma, \tag{11}$$

where $\eta = 1, \gamma = -\frac{2\delta}{D}\psi_{BI}$ for Neumann boundary conditions and $\eta = -1, \gamma = 2c_{BI}$ for Dirichlet boundary conditions. Using equation (11), the concentration $c_{IP}^n$ at the IP is represented by a function of the concentration $c_{IN}^n$ at the IN or IBN and a constant defined by the boundary condition.

The goal of the algorithm is to find the best interpolation function to estimate the value at the GIP. The interpolation function can be expressed as

$$\tilde{c}(\mathcal{P}_s) = \sum_{i=1}^{m} \mathcal{G}_i(\mathcal{P}_s)\alpha_i = \mathcal{G}^T(\mathcal{P}_s)\mathcal{A}, \tag{12}$$

where $\mathcal{P}_s$ is the coordinate of the point at $(x_s, y_s)$, $\tilde{c}(\mathcal{P}_s)$ is the interpolated value of concentration $c(\mathcal{P}_s)$ at the point $\mathcal{P}_s$, $\mathcal{g}_i$ and $\alpha_i$ are $m$ basic functions and associated coefficients, respectively. The vector forms are expressed as

$$\mathcal{G}(\mathcal{P}) = [\mathcal{g}_1(\mathcal{P}), \mathcal{g}_2(\mathcal{P}), \ldots, \mathcal{g}_m(\mathcal{P})]^T \tag{13}$$

and

$$\mathcal{A} = [\alpha_1, \alpha_2, \ldots, \alpha_m]^T. \tag{14}$$

For two-dimensional problems, the basis functions can be selected from but not limited to:

Linear: $\quad \mathcal{G}(x,y) = [1, x, y]^T \tag{15}$

Bilinear: $\quad \mathcal{G}(x,y) = [1, x, y, xy]^T \tag{16}$

Quadratic: $\quad \mathcal{G}(x,y) = [1, x, y, x^2, xy, y^2]^T \tag{17}$

Incomplete Quartic: $\quad \mathcal{G}(x,y) = [1, x, x^2, y, xy, x^2y, y^2, xy^2, x^2y^2]^T \tag{18}$

Cubic: $\quad \mathcal{G}(x,y) = [1, x, y, x^2, xy, y^2, x^3, x^2y, xy^2, y^3]^T \tag{19}$

Quartic: $\quad \mathcal{G}(x,y) = [1, x, y, x^2, xy, y^2, x^3, x^2y, xy^2, y^3, x^4, x^3y, x^2y^2, xy^3, y^4]^T \tag{20}$

Bicubic: $\quad \mathcal{G}(x,y) = [1, x, x^2, x^3, y, xy, x^2y, x^3y, y^2, xy^2, x^2y^2, x^3y^2, y^3, xy^3, x^2y^3, x^3y^3]^T \tag{21}$

In general, higher dimensional basis functions can approximate the concentration distribution near the boundary better than lower dimensional basis functions. However, higher dimensional basis functions also require more data points to get accurate and stable solutions than lower dimensional basis functions. The basis functions (15) to (21) are compared in the numerical example below.

The coefficients $\mathcal{A}$ in equation (12) can be obtained by minimizing the errors of the interpolated values at known points in a least-squares sense. The total interpolation error can be expressed as

$$\mathcal{E}_{INT}(\mathcal{A}) = \sum_{s=1}^{N} \mathcal{W}(\mathcal{P}_{GIP} - \mathcal{P}_s)[\tilde{c}(\mathcal{P}_s) - c(\mathcal{P}_s)]^2 = \sum_{s=1}^{N} \mathcal{W}(\mathcal{P}_{GIP} - \mathcal{P}_s)[\mathcal{G}^T(\mathcal{P}_s)\mathcal{A} - c(\mathcal{P}_s)]^2, \tag{22}$$

where $\mathcal{P}_s$ is the coordinate of each point used for interpolation, and $\mathcal{W}(\mathcal{P}_{GIP} - \mathcal{P}_s)$ is the weight function depending on the distance between GBI and point $s$. $\tilde{c}(\mathcal{P}_s) - c(\mathcal{P}_s)$ represents the interpolation error at point $s$. In the MLS and CMLS algorithms [27] [31], only the INs and IBNs in $\Omega_s$ are included in equation (22). In the proposed ECMLS algorithm, the associated IPs are included with INs and IBNs in equation (22). Furthermore, the associated BIs are also included if Dirichlet boundary conditions are used. The weight function assigns more weights for the points closer to the GIP in the evaluation of the total interpolation error. In the numerical implementation, a smooth and monotonically decreasing cubic spline function [31] [32] can be used as the weight function.

$$\mathcal{W}(\mathcal{P}_{GIP} - \mathcal{P}_s) = \begin{cases} 1 - 6\left(\frac{r_s}{\beta}\right)^2 + 6\left(\frac{r_s}{\beta}\right)^3, & if \ \frac{r_s}{\beta} \leq 0.5 \\ 2 - 6\left(\frac{r_s}{\beta}\right) + 6\left(\frac{r_s}{\beta}\right)^2 - 2\left(\frac{r_s}{\beta}\right)^3, & if \ 0.5 < \frac{r_s}{\beta} \leq 1.0 \\ 0, & if \ 1.0 < \frac{r_s}{\beta} \end{cases} \quad (23)$$

where $r_s$ is the normalized radius representing the distance between point $s$ and GIP

$$r_s = \frac{\|\mathcal{P}_{GIP} - \mathcal{P}_s\|}{\sqrt{\Delta x^2 + \Delta y^2}}, \quad (24)$$

and $\beta$ is a scaling parameter. The weight is zero for points $\beta\sqrt{\Delta x^2 + \Delta y^2}$ or larger distance away from the GIP (see Fig. 1 for implication of $\beta$). On the other hand, a cosine weight function has also been used [27] [30]

$$\mathcal{W}(\mathcal{P}_{GIP} - \mathcal{P}_s) = \begin{cases} \frac{1}{2}\left[1 + \cos\left(\pi \cdot \frac{r_s}{\beta}\right)\right], & if \ 0 \leq \frac{r_s}{\beta} \leq 1.0 \\ 0, & if \ 1.0 < \frac{r_s}{\beta} \end{cases} \quad (25)$$

Alternatively, the weight function can be chosen to be inversely proportional to a power of the distance [26]

$$\mathcal{W}(\mathcal{P}_{GIP} - \mathcal{P}_s) = r_s^p, \quad (26)$$

where $p$ is typically a negative number. The inverse distance weight function is most frequently used with $p = -1$.

The scaling parameter $\beta$ stretches the cosine and cubic spline weight functions laterally (Fig. 2(a)). Both weight functions are zero when the normalized radius $r_s > \beta$. For the same value of $\beta$, the cubic spline weight function (23) is steeper than the cosine weight function (25) until $r_s$ approaches $\beta$. The power of distance weight functions (Fig. 2(b)) are much steeper than the cosine and cubic spline weight functions when $r_s$ is small. For larger $r_s$, the power of distance weight functions provide low weight values and extend to large $r_s$ regions without bounds. More negative power $p$ makes the power of distance weight function (26) steeper in small $r_s$ regions and lower weight in large $r_s$ regions.

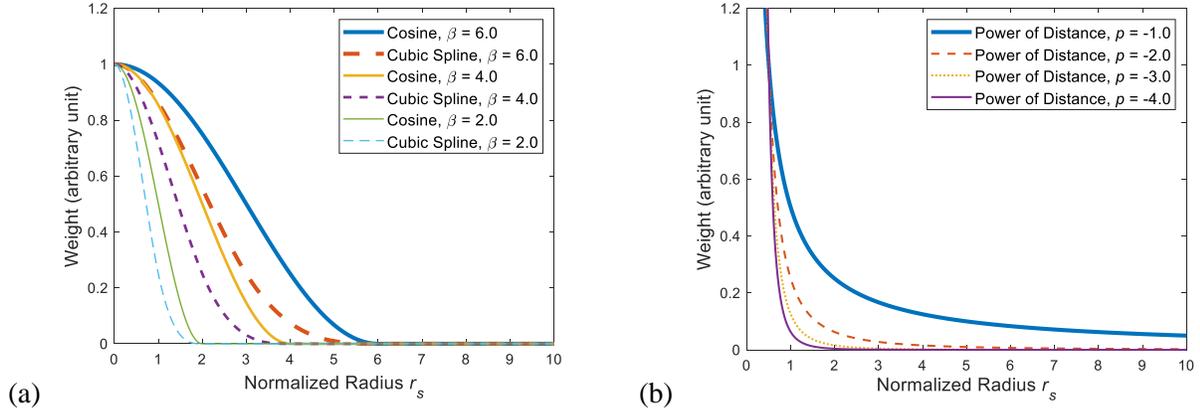

Fig. 2. Weight function comparison: (a) cosine and cubic spline weight functions; (b) power of distance weight functions normalized to 1 at $r_s = 0.5$.

Minimizing the total interpolation error $\mathcal{E}_{INT}(\mathcal{A})$ in equation (22) using the MLS algorithm can have matrix regularity issues [27]. The CMLS algorithm introduces a constraint of the boundary condition at GIP [31].

$$\mathcal{E}_{GIP}(\mathcal{A}) = [\tilde{c}(\mathcal{P}_{GIP}) - c(\mathcal{P}_{GIP})]^2. \tag{27}$$

Using equation (11) to approximate $c(\mathcal{P}_{GIP})$, equation (27) can be expressed as

$$\mathcal{E}_{GIP}(\mathcal{A}) \approx [\tilde{c}(\mathcal{P}_{GIP}) - (\eta \tilde{c}(\mathcal{P}_{GP}) + \gamma)]^2 = [\mathcal{G}^T(\mathcal{P}_{GIP})\mathcal{A} - \eta \mathcal{G}^T(\mathcal{P}_{GP})\mathcal{A} - \gamma]^2, \tag{28}$$

The total error $\mathcal{E}_{TOT}(\mathcal{A})$ includes the interpolation error $\mathcal{E}_{INT}(\mathcal{A})$ and the boundary constraint $\mathcal{E}_{GIP}(\mathcal{A})$.

$$\mathcal{E}_{TOT}(\mathcal{A}) = \mathcal{E}_{INT}(\mathcal{A}) + \kappa \mathcal{E}_{GIP}(\mathcal{A}), \tag{29}$$

where $\kappa$ is a preassigned penalty parameter which affects the strength of the constraint (27). The penalty parameter $\kappa$ is adjusted to improve the matrix regularity in solving for the interpolation coefficients.

Minimizing the total error function $\mathcal{E}_{TOT}(\mathcal{A})$ of equation (29) with respect to $\mathcal{A}$ results in the following system of equations:

$$\mathbb{M}\mathcal{A} = \mathcal{B}, \tag{30}$$

where

$$\mathbb{M} = \mathbb{G}^T \mathbb{W} \mathbb{G} + \kappa \mathcal{D} \mathcal{D}^T, \qquad \mathcal{B} = \mathbb{G}^T \mathbb{W} \mathcal{C} + \kappa \gamma \mathcal{D}. \tag{31}$$

In equation (31), $\mathbb{G}$ is the matrix of basis functions evaluated at the nodes located at $\mathcal{P}_s$, $\mathbb{W}$ is the diagonal matrix formed by weight functions $\mathcal{W}(\mathcal{P}_{GIP} - \mathcal{P}_s)$, $\mathcal{D}$ is the vector of boundary condition mismatch, and $\mathcal{C}$ is the vector of concentrations $c(\mathcal{P}_s)$.

$$\mathbb{G} = \begin{bmatrix} \mathcal{G}^T(\mathcal{P}_1) \\ \mathcal{G}^T(\mathcal{P}_2) \\ \vdots \\ \mathcal{G}^T(\mathcal{P}_N) \end{bmatrix} = \begin{bmatrix} g_1(\mathcal{P}_1) & g_2(\mathcal{P}_1) & \cdots & g_m(\mathcal{P}_1) \\ g_1(\mathcal{P}_2) & g_2(\mathcal{P}_2) & \cdots & g_m(\mathcal{P}_2) \\ \vdots & \vdots & \vdots & \vdots \\ g_1(\mathcal{P}_N) & g_2(\mathcal{P}_N) & \cdots & g_m(\mathcal{P}_N) \end{bmatrix} \quad (32)$$

$$\mathbb{W} = \begin{bmatrix} \mathcal{W}(\mathcal{P}_{GIP} - \mathcal{P}_1) & 0 & \cdots & 0 \\ 0 & \mathcal{W}(\mathcal{P}_{GIP} - \mathcal{P}_2) & \cdots & 0 \\ \vdots & \vdots & \ddots & \vdots \\ 0 & 0 & \cdots & \mathcal{W}(\mathcal{P}_{GIP} - \mathcal{P}_N) \end{bmatrix} \quad (33)$$

$$\mathcal{D} = \begin{bmatrix} g_1(\mathcal{P}_{GIP}) - \eta g_1(\mathcal{P}_{GP}) \\ g_2(\mathcal{P}_{GIP}) - \eta g_2(\mathcal{P}_{GP}) \\ \vdots \\ g_m(\mathcal{P}_{GIP}) - \eta g_m(\mathcal{P}_{GP}) \end{bmatrix} \quad (34)$$

$$\mathcal{C} = \begin{bmatrix} c(\mathcal{P}_1) \\ c(\mathcal{P}_2) \\ \vdots \\ c(\mathcal{P}_N) \end{bmatrix} \quad (35)$$

The optimal coefficient $\mathcal{A}$ that satisfies least-squares error can be found by solving equation (30). The solution can be plugged into equation (12) to estimate the concentration at the GIP.

$$c_{GIP}^n \approx \tilde{c}(\mathcal{P}_{GIP}) = \mathcal{G}^T(\mathcal{P}_{GIP})\mathcal{A} = \mathcal{G}^T(\mathcal{P}_{GIP})\mathbb{M}^{-1}\mathcal{B} = \mathcal{G}^T(\mathcal{P}_{GIP})(\mathbb{G}^T\mathbb{W}\mathbb{G} + \kappa\mathcal{D}\mathcal{D}^T)^{-1}(\mathbb{G}^T\mathbb{W}\mathcal{C} + \kappa\gamma\mathcal{D}) \quad (36)$$

The concentration at the GP can then be found through the boundary condition (11). Using equations (36) and (11), the concentration at the GP can be represented by a linear combination of concentrations at the INs and IBNs in the interpolation stencil domain $\Omega_s$. The coefficients in equation (36) depend only on the geometry of the grid points and the boundary condition. Therefore, they can be pre-calculated. Equations (36) can be plugged into the finite difference equation (7) to evaluate the associated $c_{IBN}^{n+1}$.

In equation (36), $\mathbb{M}$ is an $m$ x $m$ matrix and $\mathbb{G}$ is an $N$ x $m$ matrix. In the MLS algorithm, if $N$ is not sufficiently large compared to $m$, $\mathbb{M}$ can become singular [27] [30]. The CMLS algorithm adds the $\kappa\mathcal{D}\mathcal{D}^T$ term to improve regularity of $\mathbb{M}$ [31]. The proposed ECMLS algorithm adds the IPs into equation (22), effectively doubling $N$ without increasing the size of domain $\Omega_s$. It is confirmed by numerical examples below that the ECMLS algorithm greatly improves the regularity of $\mathbb{M}$.

## 3. Numerical Examples

To simulate diffusion in a cryo-electron microscopy (cryo-EM) grid foil hole [33], the finite difference equation (7) in Cartesian coordinates is used to model the diffusion process with a circular impermeable boundary, i.e., the Neumann boundary condition with $\psi_{BP}(\Gamma) = 0$. This investigation focuses on the accuracy of the boundary model. The simulation model is compared against an analytical solution with a point source in the center of the circular hole in order to investigate the effectiveness of the algorithms and optimization of simulation parameters.

Fick's second law of diffusion in two-dimensional homogeneous space (1) can be expressed in polar coordinates as

$$\frac{\partial c}{\partial t}(r,\theta,t) = D\left[\frac{\partial^2 c}{\partial r^2}(r,\theta,t) + \frac{1}{r}\frac{\partial c}{\partial r}(r,\theta,t) + \frac{1}{r^2}\frac{\partial^2 c}{\partial \theta^2}(r,\theta,t)\right], \quad (37)$$

where $c(r,\theta,t)$ is the concentration in polar coordinates. If the problem has a radial symmetry, equation (37) can be simplified as

$$\frac{\partial c}{\partial t}(r,t) = D\left[\frac{\partial^2 c}{\partial r^2}(r,t) + \frac{1}{r}\frac{\partial c}{\partial r}(r,t)\right]. \quad (38)$$

Assume there is an impermeable boundary at $r = R$ restricting the flow of particles, the diffusion equation would have a zero-flux boundary condition

$$\frac{\partial c}{\partial t}(R,t) = 0. \quad (39)$$

With an initial condition $c(r,0) = f(r)$, the analytical solution to the partial differential equation (38) and (39) is [34]

$$c(r,t) = \frac{1}{\pi R^2}\left[\int_0^R r'f(r')dr' + \sum_{n=1}^{\infty} \exp(-D\alpha_n^2 t)\frac{J_0(r\alpha_n)}{J_0^2(R\alpha_n)}\int_0^R r'f(r')J_0(r'\alpha_n)dr'\right], \quad (40)$$

where $\alpha_n$ are the roots of $J_1(R\alpha_n) = 0$, and $J_0$ and $J_1$ are Bessel functions of the first kind. With an initial condition $c(r,0) = \delta(r)$, where $\delta(r)$ is the Dirac delta distribution at $r = 0$, the analytical solution (40) can be simplified as

$$c_{Ana\_BC}(r,t) = \frac{1}{\pi R^2}\left[1 + \sum_{n=1}^{\infty} \exp(-D\alpha_n^2 t)\frac{J_0(r\alpha_n)}{J_0^2(R\alpha_n)}\right]. \quad (41)$$

Furthermore, with an initial condition $c(r,0) = \delta(r)$, and on an infinite plane surface, the analytical solution to the concentration distribution is [34]

$$c_{Ana\_NoBC}(r,t) = \frac{1}{4\pi Dt}\exp\left(-\frac{r^2}{4Dt}\right). \quad (42)$$

These analytical solutions were used for comparisons with the numerical simulation. The numerical simulation parameters are listed in Table 1. Three Cartesian grid sizes – 100x100, 200x200, and 400x400 – were used in the investigation with the corresponding $\Delta x$, $\Delta y$, and $\Delta t$ listed in Table 1. The diameter of the cryo-EM grid foil hole is 1 µm, and the diffusion coefficient is $10^{-10}$ m²/sec (typical value for small molecules such as ATP in water). The corresponding diffusion time $R^2/4D$ for the radius, 0.5 µm, is 625 µsec.

Table 1 Numerical simulation parameters.

| Parameter | Symbol | Value | Unit |
| --- | --- | --- | --- |
| Grid foil hole radius | $R$ | 0.5 | μm |
| Diffusion coefficient | $D$ | 1e-10 | m$^2$/sec |
| Grid size | | 100x100, 200x200, 400x400 | |
| Spatial step size | $\Delta x, \Delta y$ | 0.01, 0.005, 0.0025 | μm |
| Time step size | $\Delta t$ | 0.1, 0.025, 0.00625 | μsec |

The diffused concentration distribution with the initial condition $c(r, 0) = \delta(r)$ has a peak at the center and zero elsewhere at $t = 0$. The peak at the center spreads out radially over time. The diffusion length given by

$$L_D = \sqrt{4Dt} \tag{43}$$

characterizes the spreading of the center peak. Fig. 3 shows analytically calculated and numerically simulated concentration distribution at different time points. At $t = 1$ μsec, the diffusion length $L_D = 0.02$ μm roughly represent the half-width of the concentration peak. The diffusion length increases at the rate of square root of time.

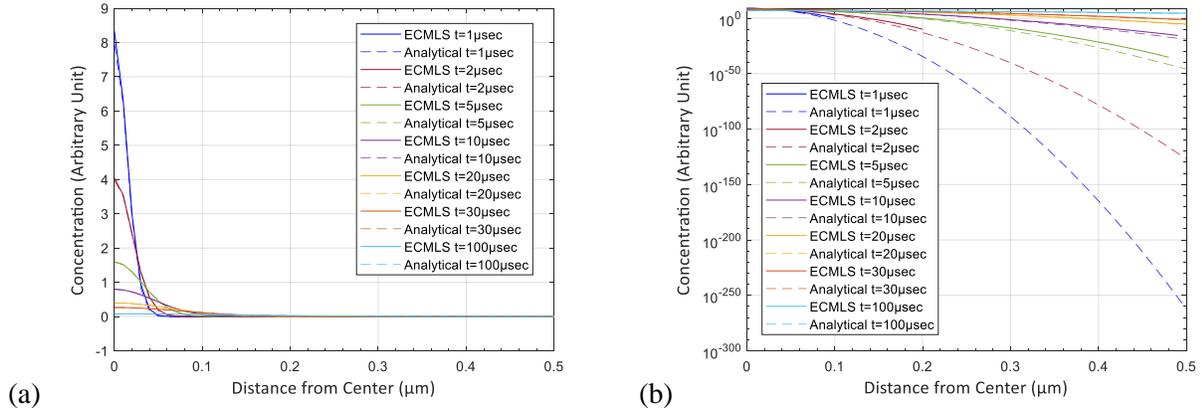

Fig. 3. Diffused concentration distribution versus distance from center over time using Cartesian coordinates with 100x100 grid: (a) linear scale; (b) logarithmic scale.

At locations several $L_D$ away from the center, the concentration values from equation (41) are very small non-zero numbers at $t > 0$. As shown in Fig. 3, these concentration values, even though small in absolute terms, can grow very rapidly in relative terms. Fig. 4 shows calculated concentration values at the boundary $x = R$ and three neighboring positions near the boundary ($R - \Delta x, R - 2\Delta x, R - 3\Delta x$) using equation (41) along the $x$ axis of a 100x100 Cartesian grid. The concentration can grow several orders of magnitude per μsec in the first few μsec. The concentration near the boundary still grows at about 10 times per μsec at the 15 μsec mark. This rapid growth is a challenge for the finite difference equation (7) to model. Fig. 4 also shows the relative concentration step values defined by

$$CS_{rel} = \frac{c_{NeighborPoint} - c_{BoundaryPoint}}{c_{BoundaryPoint}}. \tag{44}$$

The relative concentration step values at the neighbors of the boundary points can be very large in the first few μsec. The second neighbor step can still be greater than 10 at the 15 μsec mark. These large relative concentration steps cause large errors when equation (12) is used to interpolate for the concentration at the GIP. The large initial concentration growth rate is the nature of the diffusion process. On the other hand, the large relative concentration step values at the boundary can be reduced by increasing grid size. Fig. 5 shows calculated concentration values near the boundary using equation (41) with 400x400 Cartesian grid. The relative concentration step values with 400x400 grid are much lower than those with 100x100 grid. It is worth noting that the Multiprecision Computing Toolbox by Advanpix (Yokohama, Japan) had to be used to evaluate equation (41) with sufficient accuracy to calculate the concentration at $10^{-260}$ level where some conventional libraries have difficulties evaluating equation (41) below $10^{-10}$ level.

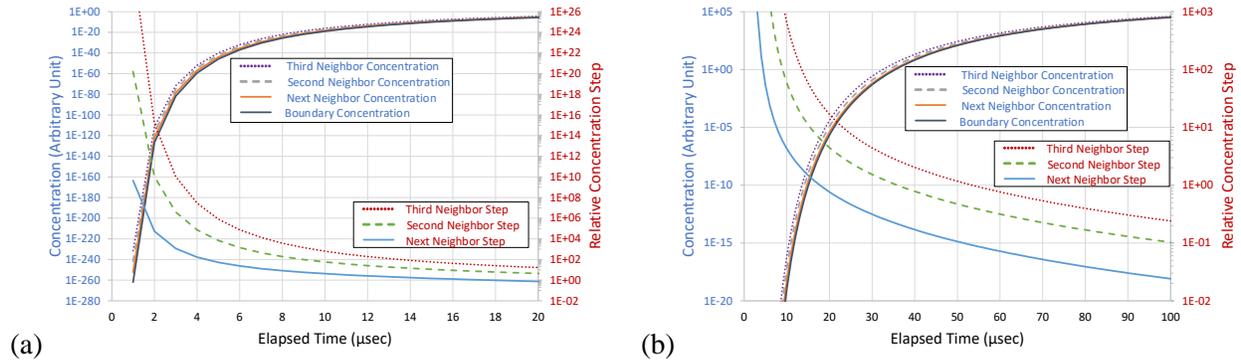

Fig. 4. Analytically calculated concentration and relative concentration step near boundaries versus time using Cartesian coordinates with 100x100 grid: (a) 1 μsec to 20 μsec view; (b) 1 μsec to 100 μsec view.

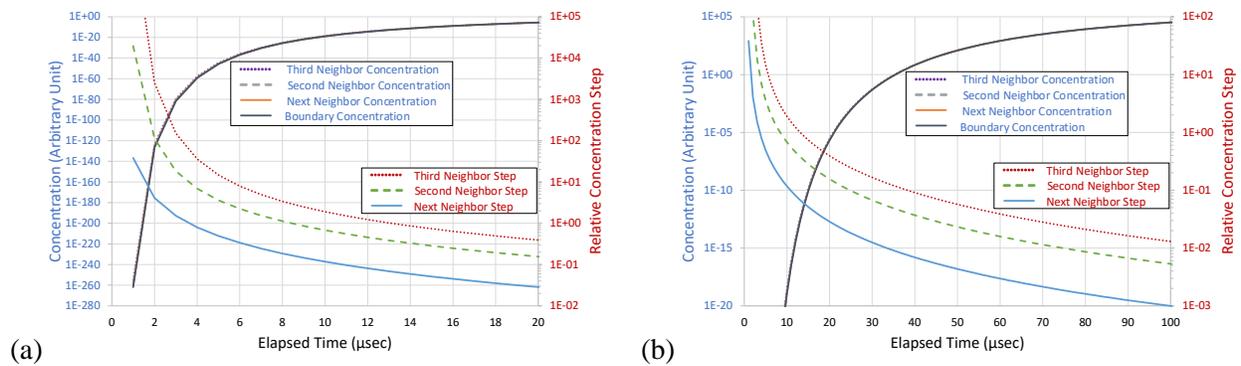

Fig. 5. Analytically calculated concentration and relative concentration step near boundaries versus time using Cartesian coordinates with 400x400 grid: (a) 1 μsec to 20 μsec view; (b) 1 μsec to 100 μsec view.

The finite difference equation (7) does not always model the diffusion equation (1) accurately. Fig. 3 also shows the numerically simulated concentration using finite difference equation (7) with the ECMLS boundary algorithm. The simulated concentration values are zero at locations more than 0.1 μm away from the center at $t = 1$ μsec. This is unrelated to the effect of diffusion length. Due to the discrete nature of finite difference method, the zero initial values at the non-center points can only be affected by the impulse release at $t = 0$ at a propagation rate of one $\Delta x$ per $\Delta t$ using equation (7). Therefore, the finite difference modeled concentration values for points outside of a radius of $t \cdot \Delta x/\Delta t$ will remain zero until the "wavefront" hits that location. In the other words, there is a null period where the analytically calculated concentration is non-zero but the finite difference modeled concentration is zero. The null period is given by

$$t_{null} = \frac{d}{\Delta x/\Delta t}, \tag{45}$$

where $d$ is the distance from center. This is unrelated to the diffusion time $d^2/4D$.

Fig. 6 shows the analytically calculated concentration at the boundary and the absolute errors and relative errors of the finite difference models. With the parameters specified in Table 1, it takes about 5 μsec for the finite difference modeled wavefront of the particles to reach the boundary in the 100x100 grid case (Fig. 6). During this null period, the relative concentration error

$$CError = \left| \frac{c_{FD\_BC} - c_{Ana\_BC}}{c_{Ana\_BC}} \right| \tag{46}$$

is 100% since the finite difference method calculated concentration $c_{FD\_BC}$ is zero, where $c_{Ana\_BC}$ is the analytically calculated concentration using equation (41). The duration of this null period can be reduced by reducing $\Delta t$.

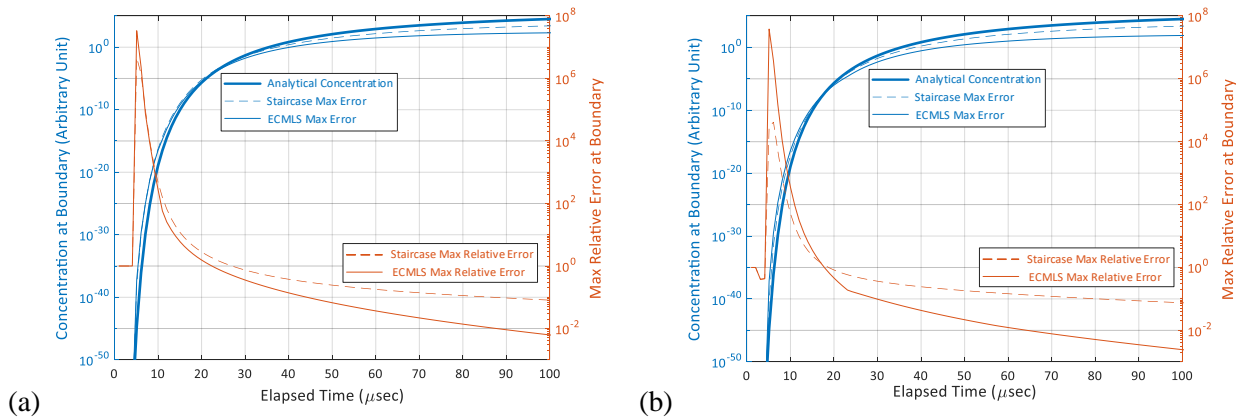

Fig. 6. Max absolute error and max relative error versus time at the boundary for a finite difference diffusion model with a point release in Cartesian coordinates with 100x100 grid: (a) uncompensated errors; (b) errors compensated for inherent finite difference induced errors. ECMLS parameters: incomplete quartic basis, cubic spline weight function, $\beta = 2.75$, $\kappa = 100$.

Immediately after the null period, the finite difference method will generate a non-zero concentration at the boundary, but this concentration value $c_{FD\_BC}$ is much higher than the analytically calculated $c_{Ana\_BC}$, resulting in the large spike of relative concentration error right after 5 μsec elapsed time in Fig. 6(a). The relative error continues to drop due to reduced error resulted from reduced relative concentration step.

The errors shown in Fig. 6(a) consists of the errors associated with the imperfect boundary model as well as the errors inherent to the finite difference model (7). Fig. 7 shows snapshots of the errors of the finite difference model (7) with ECMLS boundary model at $t = 30$ μsec. The uncompensated relative errors (Fig. 7(a)) using equation (46) without taking absolute values show positive and negative error patterns in large areas where the concentration is still orders of magnitude lower than the peak at the center. These errors are mostly inherent to the finite different equation (7) and are mostly not caused by the boundary model. To evaluate the effectiveness of boundary models, the inherent finite difference model induced errors are compensated for by using

$$CError_{comp} = \left|\frac{(c_{FD\_BC} - c_{Ana\_BC}) - (c_{FD\_NoBC} - c_{Ana\_NoBC})}{c_{Ana\_BC}}\right|, \quad (47)$$

where $c_{FD\_NoBC}$ is the finite difference calculated concentration using a grid size that is large enough to emulate a reflection-free model, and $c_{Ana\_NoBC}$ is the boundary-free analytically calculated concentration using equation (42). The compensated errors (Fig. 7(b)) using equation (47) without taking absolute values show that the errors caused by the boundary model are confined near the vicinity of the boundaries.

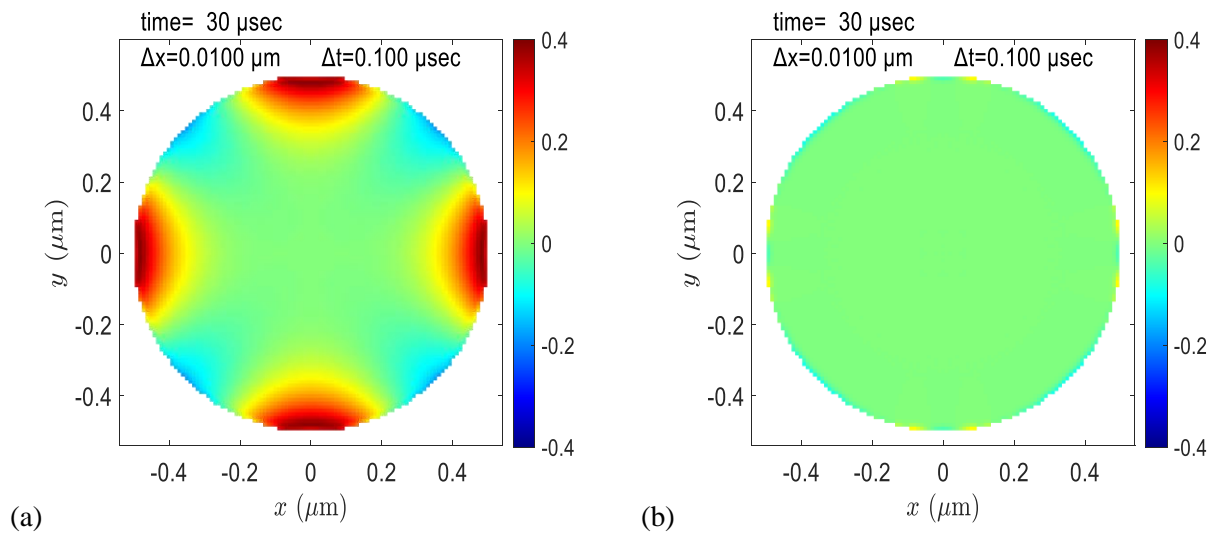

Fig. 7. Relative error of finite difference diffusion model in Cartesian coordinates with 100x100 grid at 30 μsec: (a) uncompensated errors; (b) errors compensated for inherent finite difference induced errors. ECMLS parameters: incomplete quartic basis, cubic spline weight function, $\beta = 2.75$, $\kappa = 100$.

In Fig. 6(b) the spike of compensated errors immediately after the null period is almost as high as the one of uncompensated errors for the ECMLS algorithm, even with optimized parameters in Fig. 6. This means the ECMLS algorithm produces more errors than the inherent finite difference produced errors in this region up to about 15 μsec where the relative concentration steps are very large. On the other hand, the spike of compensated errors is much smaller than that of uncompensated errors for the plain staircase model. This means the staircase model which uses only the IBN and one neighbor IN to model the boundary produces lower error than the inherent finite difference produced errors in this region. The compensated errors for the ECMLS algorithm are lower than the staircase model from 15 μsec onwards (Fig. 6(b)). This means the ECMLS algorithm which uses more neighbor points for interpolation performs better than the staircase model in the region where the relative concentration steps are lower. Fig. 6(b) also shows that the compensated error at the boundary grows higher as the diffused concentration value becomes higher. However, the compensated max error does not grow as fast as the diffused concentration value. Therefore, the compensated error (47) relative to the total concentration at the boundary becomes lower as time increases. In the region where $t > 30$ μsec, the ECMLS algorithm produces much smaller errors than the inherent finite difference induced errors. In this region, the compensated relative errors in Fig. 6(b) is much lower than the uncompensated relative errors in Fig. 6(a).

For the applications where accuracy at low concentration and early period of diffusion process is important, larger grid size can be used to reduce the relative concentration step to improve accuracy. Fig. 8 shows the uncompensated and compensated errors at the boundary with 400x400 grid size. The general behavior is very similar to the 100x100 grid size case shown in Fig. 6, except that the relative error at the boundary with 400x400 grid is much lower than that with 100x100 grid. For example, the peak relative error at the boundary is about 0.675x at $t = 10$ μsec with 400x400 grid where the error is about 278x with 100x100 grid at the same elapsed time. The larger spike in relative error in Fig. 8 is the result of shorter null time enabling the plotting of relative error in earlier elapsed time when the concentration at the boundary is much lower. The peak relative error spike with 100x100 grid in Fig. 6(a) is $3.36 \cdot 10^7$ at $t = 5$ μsec where the peak relative error with 400x400 grid in Fig. 8(a) is $4.18 \cdot 10^2$ at $t = 5$ μsec.

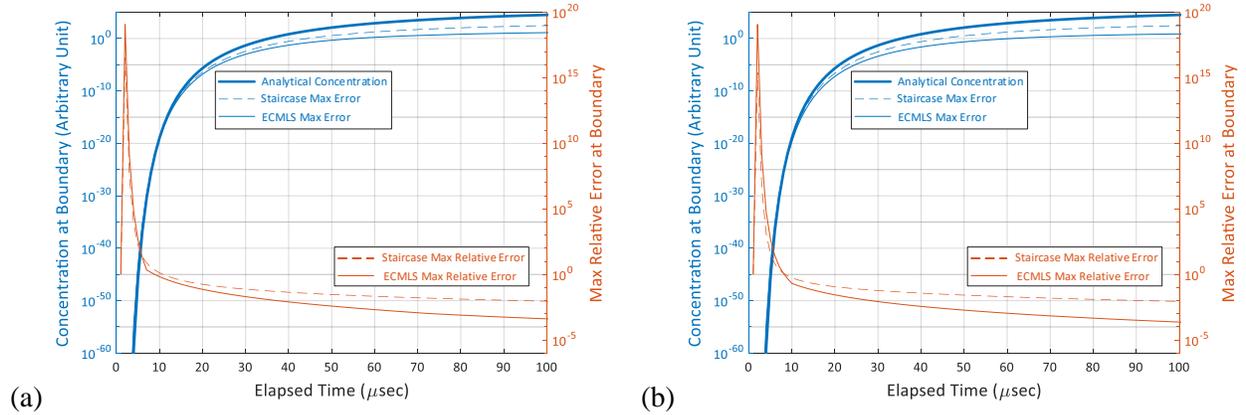

Fig. 8. Max absolute error and max relative error versus time at the boundary for a finite difference diffusion model with a point release in Cartesian coordinates with 400x400 grid: (a) uncompensated errors; (b) errors compensated for inherent finite difference induced errors. ECMLS parameters: incomplete quartic basis, cubic spline weight function, $\beta = 4.625$, $\kappa = 100$.

The compensated max relative error in Fig. 6(b) is not constant over time. However, algorithm and parameter optimization can still be performed at a fixed time. In the subsequent analysis, errors at $t = 30$ μsec are analyzed and compared in order to optimize the boundary condition model.

Using the compensated max relative error at the boundary at $t = 30$ μsec as the performance metric, various boundary models can be compared and parameters can be optimized. Fig. 9 shows the peak compensated relative error and average compensated relative error at the boundary with grid sizes 100x100, 200x200, and 400x400 at $t = 30$ μsec. Various interpolation basis functions, cubic spline weight function parameter $\beta$, and CMLS and ECMLS algorithms are compared to the staircase model. All data in Fig. 9 are generated with cubic spline weight function (23) and the penalty parameter $\kappa = 100$ in equation (29). Circular boundaries which have constant radius of curvature are used in this investigation. The use of different grid sizes affects the radius of curvature of the boundaries. In all cases, linear and bilinear basis functions do not perform better than the staircase model. As $\beta$ becomes bigger to include more points, the errors with linear and bilinear basis functions grow bigger because more points in a wider area require more complex interpolation functions to approximate. The quartic basis function does not yield better results than the rest of the basis functions for both CMLS and ECMLS algorithms, even though the results are generally better than the staircase model. The quartic basis function also causes more regularity challenges for the CMLS algorithm especially with small $\beta$ values. The cubic and bicubic basis functions do not result in good performance for the CMLS algorithm. However, the ECMLS algorithm can have good performance with both cubic and bicubic basis functions. The CMLS algorithm produces good performance only with quadratic and incomplete quartic basis functions where the incomplete quartic basis function produces slightly better results but with some regularity issue with small $\beta$ values and small radius of curvature. The ECMLS algorithm performs well with quadratic, incomplete quartic, cubic, and bicubic basis functions with no regularity issues. Among them, the incomplete quartic basis function offers the best overall performance with the ECMLS algorithm across grid sizes. The optimal $\beta$ values are 2.75, 3.625, and 4.625 for grid sizes of 100x100, 200x200, and 400x400, respectively. This means that it is more optimal to use smaller regions for interpolation when the radius of curvature of the local boundary is small.

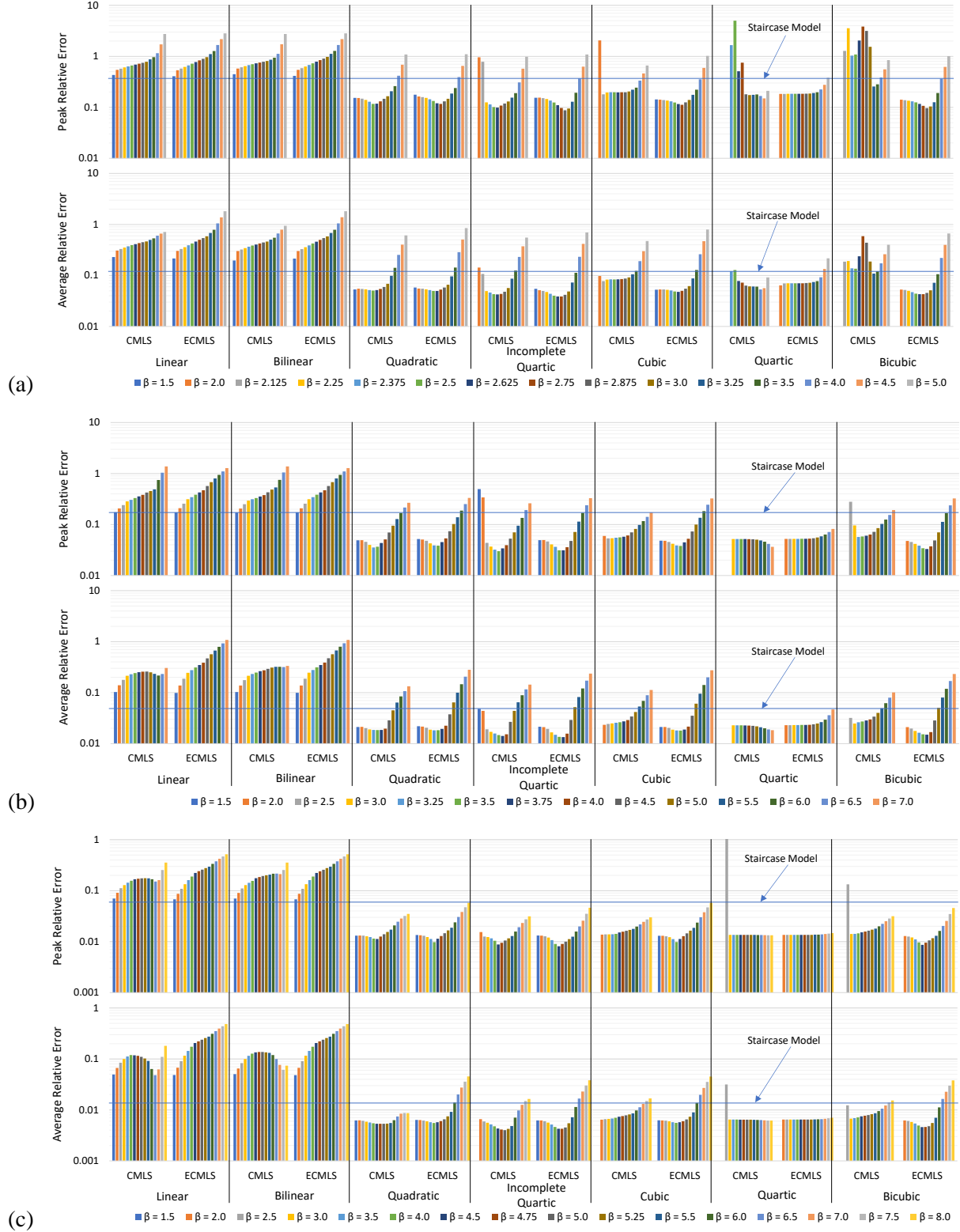

Fig. 9. Peak compensated relative error (top) and average compensated relative error (bottom) at IBN points within a circular impermeable boundary at 30 μsec elapsed time: (a) 100x100 Cartesian grid; (b) 200x200 Cartesian grid; (c) 400x400 Cartesian grid. Other parameters: cubic spline weight function with various $\beta$ values; penalty parameter $\kappa = 100$; elapsed time $t = 30$ μsec.

The penalty parameter $\kappa$ is introduced by the CMLS algorithm [31] to improve the matrix regularity of the solutions to equation (36). The proposed ECMLS algorithm also incorporated this penalty parameter. Fig. 10 shows the sensitivity of compensated relative error (47) to the penalty parameter. With $\kappa = 0$, the CMLS algorithm becomes the MLS algorithm [27]. Fig. 10 shows that the matrix regularity issue is prevalent for $\kappa < 10$ with CMLS and MLS algorithms. It takes a penalty parameter $\kappa \geq 100$ to improve matrix regularity in CMLS. However, the ECMLS algorithm does not require the penalty parameter to stabilize matrix regularity since there is no observable performance degradation at $\kappa = 0$. However, it shows slight performance improvement using larger grid size, implying larger boundary radius of curvature, with $\kappa \geq 1000$.

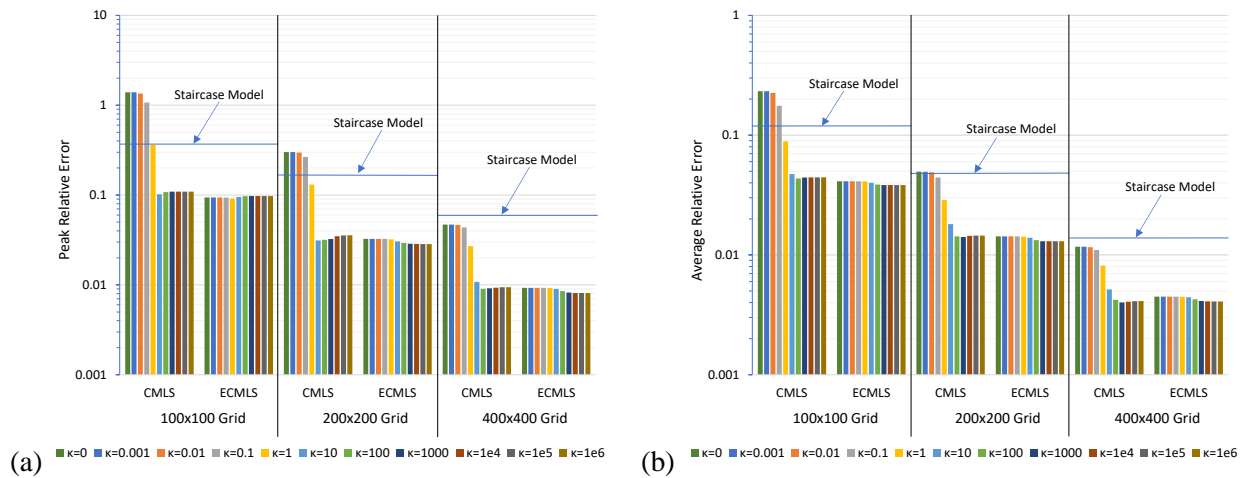

Fig. 10. Sensitivity of compensated relative error to the penalty parameter $\kappa$: (a) peak compensated relative error; (b) average compensated relative error. Other parameters: incomplete quartic basis function; cubic spline weight function with $\beta = 2.75$ (100x100), $\beta = 3.625$ (200x200), and $\beta = 4.625$ (400x400); elapsed time $t = 30$ μsec.

The effectiveness comparison and parameter optimization of the weight functions are also performed. Fig. 11, Fig. 12, and Fig. 13 show the sensitivity of compensated relative errors to the parameters of the cubic spline weight function (23), cosine weight function (25), and power of distance weight function (26), respectively. These weight function parameters effectively define the region of points used in the interpolation process. As shown in Fig. 2(a), the shape of the cosine weight function is similar to and slightly wider than that of the cubic spline weight function with the same scaling parameter $\beta$. Fig. 11 and Fig. 12 show that the performance of the cubic spline weight function is similar to the cosine weight function. The optimal $\beta$ is slightly smaller with the cosine weight function than the cubic spline weight function. This is consistent with the fact that the width of the cosine weight curve is slightly wider than that of the cubic spline weight curve. The optimal $\beta$ in both weight functions is smaller with smaller grid size implying smaller boundary radius of curvature. When $\beta$ is too small, there are not enough points in the smaller area to provide information for interpolation. When $\beta$ is too large, the concentration distribution in the larger area is too complicated for the interpolation function (12) to approximate, especially with smaller boundary radius of curvature. The ECMLS algorithm shows good numerical stability across wide range of $\beta$. However, the CMLS algorithm might have matrix regularity issues with small $\beta \leq 2$. As shown in Fig. 2(b), the shape of the power of distance weight function has a gentle roll-off with low magnitude power parameters $p$ and a very sharp drop-off with larger magnitude $p$. Fig. 13 shows that the performance for both CMLS and ECMLS algorithms with low magnitude $p$ is worse than the staircase model. With the larger 400x400 grid, it takes the $p$ value at around -3.5 to get optimal performance. With the smaller 100x100 grid, it takes the $p$ value at around -8 to get optimal performance. This result is very different from the reported choice of $p = -1$ in the other investigation [26].

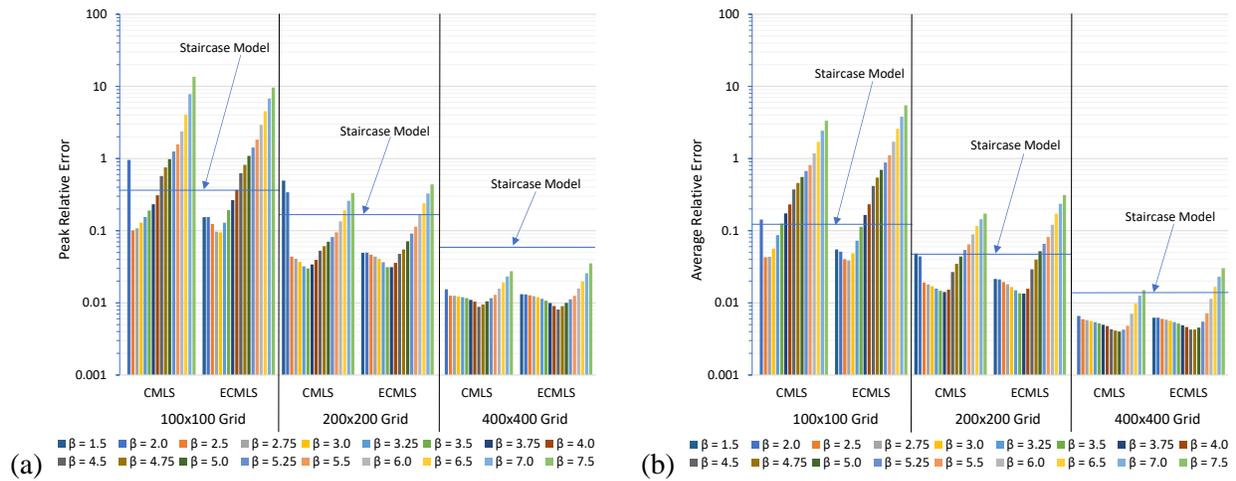

Fig. 11. Sensitivity of compensated relative error to the scaling parameter $\beta$ of cubic spline weight function: (a) peak compensated relative error; (b) average compensated relative error. Other parameters: incomplete quartic basis function; penalty parameter $\kappa = 100$; elapsed time $t = 30$ μsec.

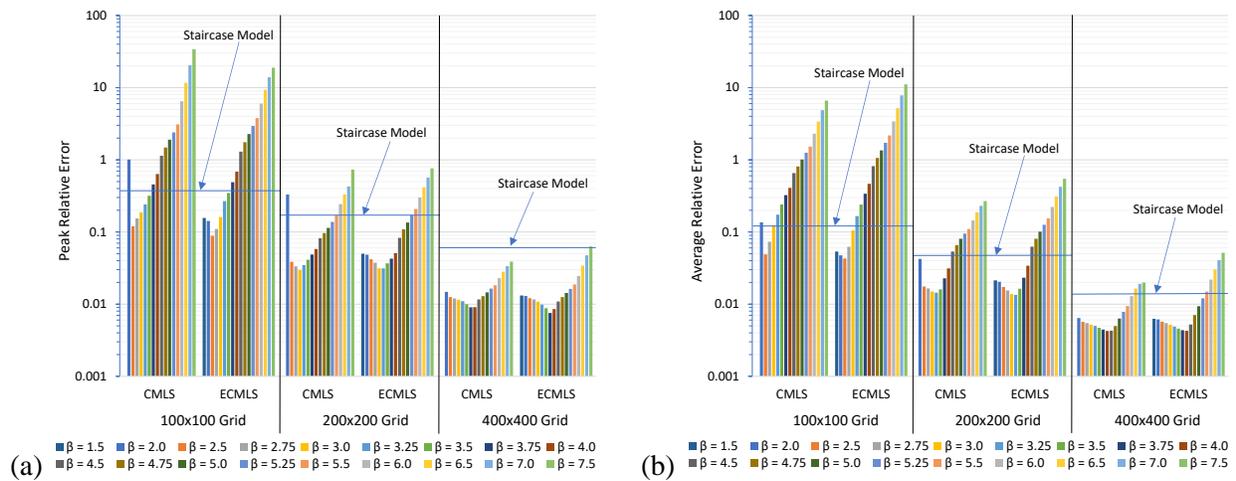

Fig. 12. Sensitivity of compensated relative error to the scaling parameter $\beta$ of cosine weight function: (a) peak compensated relative error; (b) average compensated relative error. Other parameters: incomplete quartic basis function; penalty parameter $\kappa = 100$; elapsed time $t = 30$ μsec.

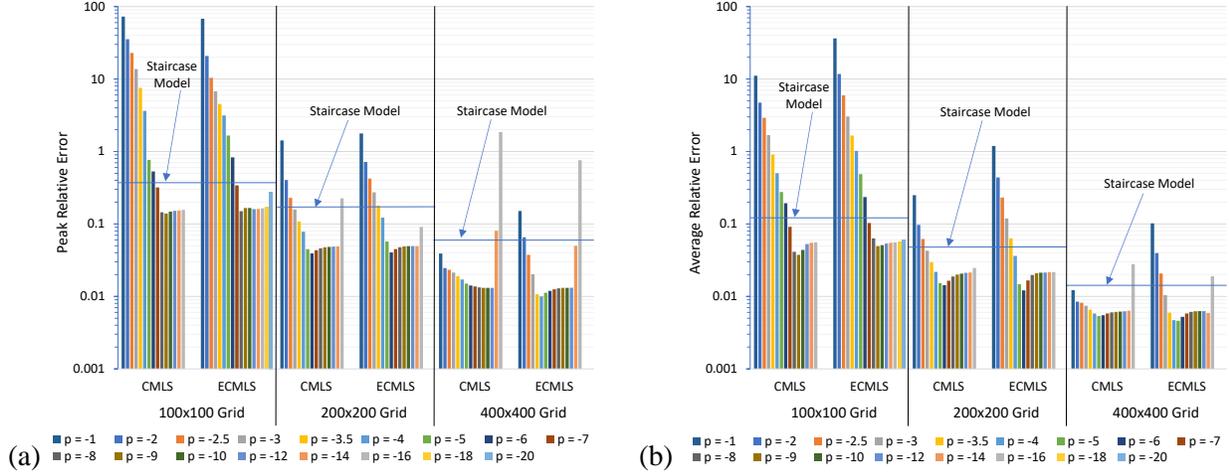

Fig. 13. Sensitivity of compensated relative error to the power parameter $p$ of power of distance weight function: (a) peak compensated relative error; (b) average compensated relative error. Other parameters: incomplete quartic basis function; penalty parameter $\kappa = 100$; elapsed time $t = 30$ μsec.

To investigate the optimization and dependency of the scaling parameter $\beta$ of cubic spline weight function, cases listed in Table 2 are simulated using incomplete quartic basis function for interpolation. A wide range of the scaling parameter $\beta$ is surveyed to find the optimal $\beta$ that minimizes the compensated relative error (47) at the boundary. Using the ECMLS algorithm, it is found that the optimal $\beta$ has strong correlation to the radius of curvature of the boundary which is the radius of the circular region. As shown in Fig. 14, the optimal $\beta$ linearly scales with the logarithm of radius of curvature and is independent of the spatial step size. Therefore, in regions with small radius of curvature at the boundary, it is better to use smaller $\beta$, implying a smaller interpolation stencil domain, to interpolate the concentration at the GIP. In regions of large radius of curvature at the boundary, using larger $\beta$ to include more data points is better. The empirically derived formula for the optimal scaling parameter $\beta$ of cubic spline weight function can be expressed as

$$\beta_{opt} = 3.2107 \cdot \log_{10}\left(\frac{R_c}{\Delta x}\right) - 2.7501, \tag{48}$$

where $R_c$ is the radius of curvature of the boundary, i.e., the radius of the cryo-EM grid hole. For applications with complex boundary structure, the scaling parameter $\beta$ of cubic spline weight function can be adjusted according to equation (48) depending on the radius of curvature of the local boundary to achieve optimal performance.

Table 2 Simulation cases used to investigate optimal scaling parameter $\beta$ of cubic spline weight function.

| Radius R (μm) | Grid Size | Spatial Step Size $\Delta x, \Delta y$ (μm) | Time Step Size $\Delta t$ (μsec) |
|---|---|---|---|
| 0.0625 | 50x50 | 0.0025 | 0.00625 |
| 0.125 | 50x50 | 0.005 | 0.025 |
| 0.25 | 50x50 | 0.01 | 0.1 |
| 0.125 | 100x100 | 0.0025 | 0.00625 |
| 0.25 | 100x100 | 0.005 | 0.025 |
| 0.5 | 100x100 | 0.01 | 0.1 |
| 0.25 | 200x200 | 0.0025 | 0.00625 |
| 0.5 | 200x200 | 0.005 | 0.025 |
| 1.0 | 200x200 | 0.01 | 0.1 |
| 0.5 | 400x400 | 0.0025 | 0.00625 |
| 1.0 | 400x400 | 0.005 | 0.025 |
| 2.0 | 400x400 | 0.01 | 0.1 |

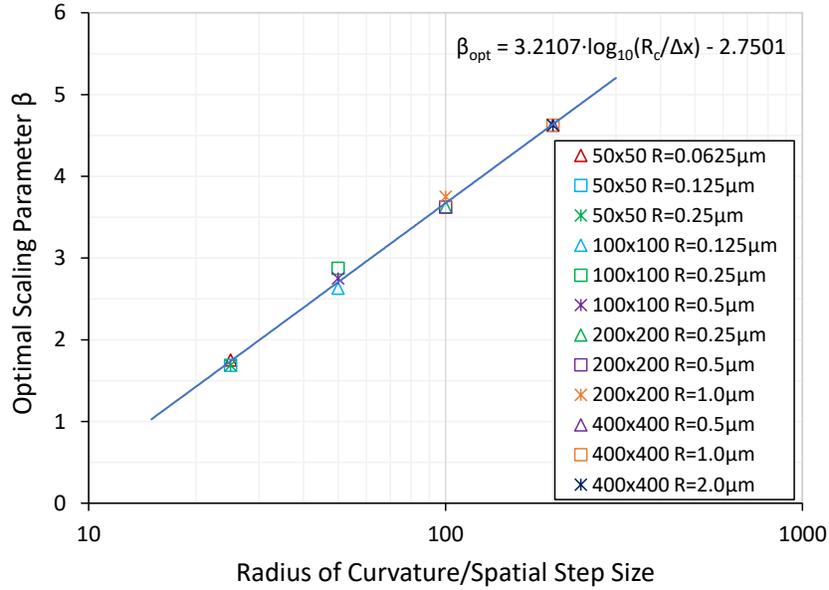

Fig. 14. Optimal scaling parameter $\beta$ of cubic spline weight function that minimizes the compensated relative error. Other parameters: ECMLS algorithm, incomplete quartic basis function; penalty parameter $\kappa = 100$.

From the results in Fig. 14, the optimal $\beta$ of cubic spline weight function scales linearly with the logarithm of radius of curvature. With the 50x50 grid where the boundary radius of curvature is $25\Delta x$, the optimal $\beta$ is around 1.6875. Fig. 15 shows the performance comparison between CMLS and ECMLS algorithms with a wide range of scaling parameter $\beta$ using a 50x50 grid. The CMLS algorithm does not produce good results with $\beta < 2.25$. However, the ECMLS algorithm does not have penalty in the small $\beta$ region. Therefore, the ECMLS algorithm could provide more optimized results than the CMLS algorithm in the region with small boundary radius of curvature.

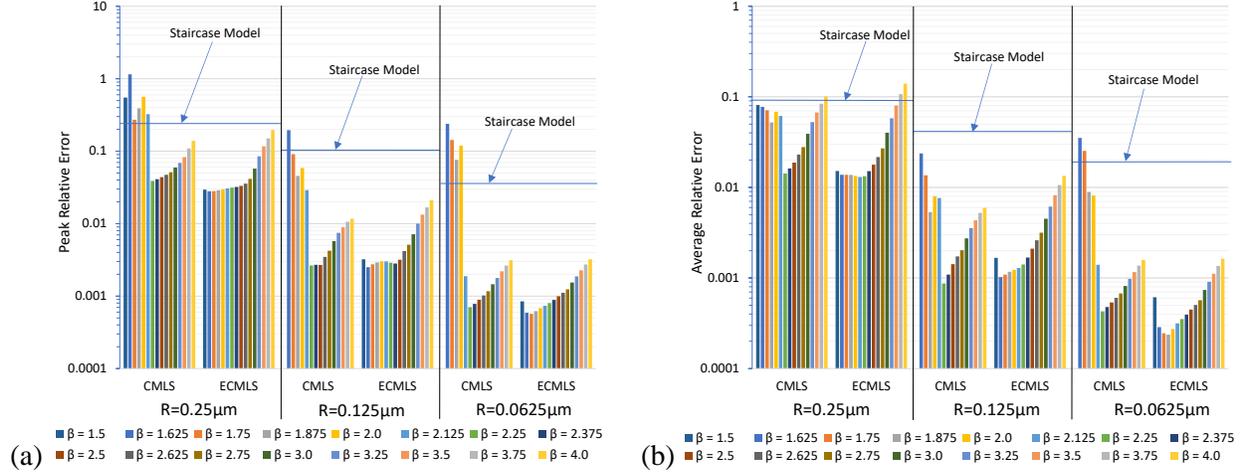

Fig. 15. Sensitivity of compensated relative error to the scaling parameter $\beta$ of cubic spline weight function: (a) peak compensated relative error; (b) average compensated relative error. Other parameters: 50x50 grid size, incomplete quartic basis function; penalty parameter $\kappa = 100$.

## 4. Conclusions

A novel sharp interface ghost-cell based immersed boundary method has been proposed and investigated with the staircase model and the constrained moving least-squares (CMLS) interpolation methods in diffusion applications. The proposed embedded constrained moving least-squares (ECMLS) algorithm incorporates the mirrored image points of all internal nodes inside the interpolation stencil domain through the associated boundary conditions to interpolate the concentration at the image point of the ghost point. The ECMLS algorithm incorporates as many boundary constraints as the number of internal points used in interpolation, where the CMLS algorithm incorporates one constraint and the MLS algorithm incorporates none. For the investigated application of substrate diffusion in grid holes with impermeable boundary in time-resolved cryo electron microscopy experiments, the ECMLS algorithm produces results better than the staircase model by using the quadratic, incomplete quartic, cubic, and bicubic basis functions. The incomplete quartic basis function yields the best performance among these basis functions. In contrast, the CMLS algorithm produces good results, comparable to the ECMLS results, only with quadratic and incomplete quartic basis functions. It is also found that both ECMLS and CMLS algorithms cannot produce results better than the staircase model when either linear or bilinear basis function is used. In the region with small boundary radius of curvature, it is found that smaller interpolation stencil domain size is more optimal in interpolating the value at the ghost image point. It is shown that the optimal radius of the interpolation stencil domain scales with the logarithm of the boundary radius of curvature. In summary, the proposed ECMLS algorithm produces lower errors at the boundary with better numerical stability over a wider range of basis functions, weight functions, boundary radius of curvature, and the penalty parameter than the CMLS and staircase algorithms.

**Acknowledgments**

This work was in part supported by SLAC National Accelerator Laboratory and Stanford University (to S.W. and P.-N. Li).
**References**

[1] D. B. Spalding, "A novel finite difference formulation for differential expressions involving both first and second derivatives," *International Journal for Numerical Methods in Engineering,* vol. 4, no. 4, pp. 551-559, July/August 1972.

[2] W. H. Press, S. A. Teukolsky, W. T. Vetterling and B. P. Flannery, "Chapter 20. Partial Differential Equations," in *Numerical Recipes*, third ed., Cambridge, Cambridge University Press, 2007, pp. 1024-1096.

[3] J. H. Ferziger, M. Perić and R. L. Street, Computational Methods for Fluid Dynamics, fourth ed., Springer, 2020.

[4] J. H. Bramble and B. E. Hubbard, "Approximation of Solutions of Mixed Boundary Value Problems for Poisson's Equation by Finite Differences," *Journal of the ACM,* vol. 12, no. 1, pp. 114-123, January 1965.

[5] A. R. Mitchell and D. F. Griffiths, The finite difference method in partial differential equations, John Wiley & Sons, 1980.

[6] K. S. Yee, "Numerical solution of initial boundary value problems involving Maxwell's equations in isotropic media," *IEEE Transactions on Antennas and Propagation,* vol. 14, no. 3, pp. 302-307, May 1966.

[7] J. Häggblad and O. Runborg, "Accuracy of staircase approximations in finite-difference methods for wave propagation," *Numerische Mathematik,* vol. 128, pp. 741-771, March 2014.

[8] R. A. Drainville, L. Curiel and S. Pichardo, "Superposition method for modelling boundaries between media in viscoelastic finite difference time domain simulations," *The Journal of the Acoustical Society of America,* vol. 146, no. 6, pp. 4382-4401, December 2019.

[9] J. F. Thompson, Z. U. A. Warsi and C. W. Mastin, "Boundary-fitted coordinate systems for numerical solution of partial differential equations—A review," *Journal of Computational Physics,* vol. 47, no. 1, pp. 1-108, July 1982.
25

[10] T.-K. Tsay, B. A. Ebersole and P. L.-F. Liu, "Numerical modelling of wave propagation using parabolic approximation with a boundary-fitted co-ordinate system," *International Journal for Numerical Methods in Engineering,* vol. 27, no. 1, pp. 37-55, July 1989.

[11] M. R. Visbal and D. V. Gaitonde, "On the Use of Higher-Order Finite-Difference Schemes on Curvilinear and Deforming Meshes," *Journal of Computational Physics,* vol. 181, no. 1, pp. 155-185, September 2002.

[12] Y. Jiang, C.-W. Shu and M. Zhang, "Free-stream preserving finite difference schemes on curvilinear meshes," *Methods and Applications of Analysis,* vol. 21, no. 1, pp. 1-30, April 2014.

[13] C. S. Peskin, "Flow patterns around heart valves: A numerical method," *Journal of Computational Physics,* vol. 10, no. 2, pp. 252-271, October 1972.

[14] D. Goldstein, R. Handler and L. Sirovich, "Modeling a No-Slip Flow Boundary with an External Force Field," *Journal of Computational Physics,* vol. 105, no. 2, pp. 354-366, April 1993.

[15] R. Mittal and G. Iaccarino, "Immersed Boundary Methods," *Annual Review of Fluid Mechanics,* vol. 37, pp. 239-262, 2005.

[16] H. S. Udaykumar, R. Mittal, P. Rampunggoon and A. Khanna, "A Sharp Interface Cartesian Grid Method for Simulating Flows with Complex Moving Boundaries," *Journal of Computational Physics,* vol. 174, pp. 345-380, 2001.

[17] J. H. Seo and R. Mittal, "A sharp-interface immersed boundary method with improved mass conservation and reduced spurious pressure oscillations," *Journal of Computational Physics,* vol. 230, pp. 7347-7363, 2011.

[18] B. Muralidharan and S. Menon, "A high-order adaptive Cartesian cut-cell method for simulation of compressible viscous flow over immersed bodies," *Journal of Computational Physics,* vol. 321, pp. 342-368, September 2016.

[19] R. P. Fedkiw, A. Marquina and B. Merriman, "An Isobaric Fix for the Overheating Problem in Multimaterial Compressible Flows," *Journal of Computational Physics,* vol. 148, pp. 545-578, 1999.

[20] Y.-H. Tseng and J. H. Ferziger, "A ghost-cell immersed boundary method for flow in complex geometry," *Journal of Computational Physics,* vol. 192, pp. 593-623, 2003.

[21] R. Ghias, R. Mittal and H. Dong, "A sharp interface immersed boundary method for compressible viscous flows," *Journal of Computational Physics,* vol. 225, no. 1, pp. 528-553, July 2007.

[22] K. Luo, Z. Zhuang, J. Fan and N. E. L. Haugen, "A ghost-cell immersed boundary method for simulations of heat transfer in compressible flows under different boundary conditions," *International Journal of Heat and Mass Transfer,* vol. 92, pp. 708-717, 2016.

[23] R. Mittal, H. Dong, M. Bozkurttas, F. Najjar, A. Vargas and A. von Loebbecke, "A versatile sharp interface immersed boundary method for incompressible flows with complex boundaries," *Journal of Computational Physics,* vol. 227, pp. 4825-4852, 2008.

[24] A. Kapahi, S. Sambasivan and H. Udaykumar, "A three-dimensional sharp interface Cartesian grid method for solving high speed multi-material impact, penetration and fragmentation problems," *Journal of Computational Physics,* vol. 241, pp. 308-332, 2013.

[25] H. Uddin, R. Kramer and C. Pantano, "A Cartesian-based embedded geometry technique with adaptive high-order finite differences for compressible flow around complex geometries," *Journal of Computational Physics,* vol. 262, pp. 379-407, 2014.

[26] C. Brehm, C. Hader and H. Fasel, "A locally stabilized immersed boundary method for the compressible Navier-Stokes equations," *Journal of Computational Physics,* vol. 295, pp. 475-504, 2015.

[27] J. H. Seo and R. Mittal, "A high-order immersed boundary method for acoustic wave scattering and low-Mach number flow-induced sound in complex geometries," *Journal of Computational Physics,* vol. 230, no. 4, pp. 1000-1019, February 2011.

[28] J. Xia, K. Luo and J. Fan, "A ghost-cell based high-order immersed boundary method for inter-phase heat transfer simulation," *International Journal of Heat and Mass Transfer,* vol. 75, pp. 302-312, 2014.

[29] J. Fernández-Fidalgo, S. Clain, L. Ramírez, I. Colominas and X. Nogueira, "Very high-order method on immersed curved domains for finite difference schemes with regular Cartesian grids," *Computer Methods in Applied Mechanics and Engineering,* vol. 360, p. 112782, March 2020.

[30] H. Luo, R. Mittal, X. Zheng, S. A. Bielamowicz, R. J. Walsh and J. K. Hahn, "An immersed-boundary method for flow–structure interaction in biological systems with application to phonation," *Journal of Computational Physics,* vol. 227, no. 22, pp. 9303-9332, November 2008.

[31] Y. Qu, R. Shi and R. C. Batra, "An immersed boundary formulation for simulating high-speed compressible viscous flows with moving solids," *Journal of Computational Physics,* vol. 354, pp. 672-691, 2018.

[32] J. Dolbow and T. Belytschko, "Numerical integration of the Galerkin weak form in meshfree methods," *Computational Mechanics,* vol. 23, pp. 219-230, April 1999.

[33] C. J. Russo and L. A. Passmore, "Progress towards an optimal specimen support for electron cryomicroscopy," *Current Opinion in Structural Biology,* vol. 37, pp. 81-89, January 2016.

[34] J. Crank, "5. Diffusion in a cylinder," in *The Mathematics of Diffusion*, second ed., Oxford University Press, 1975.